\documentclass{amsart}
\usepackage{amssymb,amsthm}
\usepackage{mathrsfs}
\def\mathcal{\mathscr}
\def\fk{\mathfrak}
\def\NN{{\mathbb N}}

\def\eps{\varepsilon}
\def\al{\alpha}

\def\del{\delta}

\def\rho{\varrho}
\def\phi{\varphi}

\newcommand{\f}[2]{\frac{#1}{#2}}
\newcommand{\be}[1]{\begin{equation}\label{#1}}
\newcommand{\ee}{\end{equation}}
\newcommand{\multsum}[2]{\sum_{{\scriptstyle #1}\atop {\scriptstyle #2}}}
\newcommand{\trisum}[3]{\sum_{{{\scriptstyle #1}\atop {\scriptstyle #2}}\atop{\scriptstyle #3}}}

\newcommand{\rf}[1]{{\rm (\ref{#1})}}

\def\cal{\mathcal}

\def\hh{{\mathbf h}}
\newcommand{\bfnu}{\mbox{\boldmath $\nu$}}
\newtheorem{lemma}{Lemma}
\newtheorem{theorem}{Theorem}
\newtheorem*{corollary}{Corollary}
\begin{document}
\title[ moderately dense sequences of primes]{Local oscillations in \\ moderately dense sequences of primes}
\author{J\"org Br\"udern and Christian Elsholtz}
\dedicatory{To Robert Tichy, on the occasion of his 60th birthday}
\address{
J\"org Br\"udern, Mathematisches Institut, Bunsenstrasse 3--5, 37073 G\"ottingen, Germany.}
\email{Joerg.Bruedern@mathematik.uni-goettingen.de}
\address{Christian Elsholtz, Institut f\"ur Analysis und Zahlentheorie, Technische Universit\"at Graz, Kopernikusgasse 24, A-8010 Graz, Austria}
\email{elsholtz@math.tugraz.at}
\thanks{The authors a grateful to CIRM at Marseille Luminy for creating a stimulating working atmosphere. The second author also likes to thank Forschungsinstitut Mathematik (FIM) at ETH Z\"urich for a very pleasant stay.}
\maketitle
\begin{abstract} The distribution of differences of consecutive members of  sequences  of primes is investigated. A quantitative measure
for oscillations among these differences is the curvature of the sequence. If the sequence is not too sparse, then sharp estimates for its curvature 
are provided.
\end{abstract}
\section{Introduction}
In an influential paper, Erd\H os and Tur\'an \cite{ET} showed that when $(p_n)$
denotes the sequence  of all prime numbers arranged in increasing order, 
then there are infinitely many sign changes among the numbers
\be{trip} p_{n+1}^2- p_n p_{n+2}.  \ee
Motivated by quantitative versions of this result due to  R\'enyi \cite{R} and
Erd\H os and  R\'enyi \cite{ER}, we develop this theme further in the context of sequences that are not too sparse. 

\begin{theorem}\label{thm0} Let  $\cal P$ be a set of primes with the property that
\be{000} \frac{(\log x)^{4/3}}{x}   \# \{p\in{\cal P}: p\le x\} \ee
tends to infinity with $x$. If $p_n$ denotes the enumeration of the set $\cal P$ in increasing order, then the sequence
\rf{trip} changes sign infinitely often.
\end{theorem}  

Our main object of study is the {\em curvature} of sequences. The idea is due to R\'enyi~\cite{R}. Consider at least three distinct points $z_1,\ldots,z_N$ in the complex plane. With the argument of a complex number 
chosen in the interval $(-\pi,\pi]$, the sum
\be{1}  \sum_{n=1}^{N-2} \Big| \arg \f{z_{n+2}-z_{n+1}}{z_{n+1}-z_{n}}\,\Big|.\ee
is referred to as the total curvature of the polygonal line connecting $z_{n-1}$ with $z_n$ for $2\le n\le N$, because this adds up the (non-negative) angles between the line segments from $z_n$ to $z_{n+1}$, and on to $z_{n+2}$. For a set of primes $\cal P$, again enumerated in increasing order as $p_n$, we take  
$z_n= n+{\mathrm i} \log p_n$, and then let 
$K_N({\cal P})$ denote the sum in \rf{1} with this special choice of $z_n$. This is the curvature of $\cal P$, truncated at $N$.

Now suppose that we knew that $K_N({\cal P})$ were unbounded. Then, if the segment 
$(\log p_n)_{n_0\le n \le N}$ is either concave or convex, then $K_N({\cal P})-K_{n_0}({\cal P})\le \f12 \pi$ which is impossible for large $N$. We conclude that the sequence 
$$ \log p_{n+2} - 2\log p_{n+1} + \log p_n$$
changes sign infinitely often, and on taking exponentials this is the same as exhibiting sign changes in the sequence \rf{trip}. In particular, Theorem \ref{thm0} will follow once we have established that $K_N({\cal P})$ is unbounded for the sequences of primes satisfying \rf{000}. Further, we see that the growth rate of $K_N({\cal P})$ is a rough measure for the oscillations in the sequence \rf{trip}. 

R\'enyi \cite{R} in 1950 considered the sequence of all primes and bounded their curvature, hereafter denoted by $K_N$, from below by
$$ K_N \gg \log\log \log N. $$
Shortly afterwards, in collaboration with Erd\H os \cite{ER} (see also \cite{P}) he determined the order of magnitude of $K_N$, now showing that
\be{2} \log N \ll K_N \ll \log N. \ee 
Their methods rely on the prime number theorem. Our concern in this paper is with estimates for the curvature that are based solely on lower bounds for the number of primes in a given sequence, such as in \rf{000}.  Before we can formulate our principal estimate, we have to set up some notation.

We work relative to an arithmetic progression. When $a,q\in\NN$ with $1\le a\le q$ and $(a,q)=1$,
let ${\cal P}_{q,a}$ denote the set of all primes $p\equiv a\bmod q$.
We refer to a subset  ${\cal P}\subset{\cal P}_{q,a} $ as {\em dense} if there are positive numbers $\delta$ and $x_0$ with the property that whenever $x\ge x_0$, then
\be{3} \# \{p\in{\cal P}: p\le x\} \ge \delta \pi(x;q,a) \ee
where as usual  $\pi(x;q,a)$ is the number of primes not exceeding $x$ in ${\cal P}_{q,a}$. 
More generally, if $\delta: [3,\infty)\to (0,1]$ is monotonically decreasing with $\delta(x)\ge (\log x)^{-1}$, and \rf{3} is satisfied with $\delta=\delta(x)$ for all $x\ge x_0$, then\footnote{It may seem unnatural to include the lower bound on $\delta$ in this definition, but more rapidly decaying functions will play no role in this paper, and it simplifies the exposition later that $\delta$ is not too small, {\em a fortiori}.} the set $\cal P$ is called $\delta$-dense (relative to $x_0$ and ${\cal P}_{q,a}$). The lower bound on $\delta$ ensures that $\cal P$ is an infinite set,  
enumerated in ascending order by $p_n$, as before. Then $K_N({\cal P})$ is defined for all $N\ge 3$.
We also put $\delta_N= \delta(p_{N})$.

\begin{theorem}\label{thm1}
Fix a number $x_0\ge 3$ and a decreasing function   $\delta: [3,\infty)\to (0,1]$ with $\delta(x)\ge (\log x)^{-1}$
for all $x\ge 3$. Then there is a sequence of natural numbers $N_0(q)$ 
with the property that for all  $N\ge N_0(q)$ and for all sets of primes  
$\cal P$ that are  $\delta$-dense relative to $x_0$ and some ${\cal P}_{q,a}$,  
one has
\be{krumm}   K_N({\cal P})\le 500 \del_{N}^{-1}\log N. \ee
If $\delta(x)^2\log x$ tends to infinity with $x$, then one also has
$$  K_N({\cal P}) \ge 10^{-8}\delta_N^3 \log N. $$
\end{theorem}

Theorem \ref{thm1} may be applied to the arithmetic progression ${\cal P}_{q,a}$ itself, with $\delta=1$. We then conclude as follows.

\begin{corollary} With $N_0(q)$ as in the preceding Theorem, for $N\ge N_0(q)$ one has 
$$ 10^{-8}  \log N  \le  K_N({\cal P}_{q,a})\le 500\log N .$$
\end{corollary}
This contains \rf{2} as a very special case. Note that here as well as in Theorem \ref{thm1} no effort has been made to optimise the numerical constants.

When $\delta$ decays it is important to have at hand a lower bound for $\delta_N$. One has 
\be{delta} \delta_N \ge \delta(4\phi(q)N(\log N)^{2}) \ee
for all large $N$.
We show this in passing, in Section 3 below. In particular,
 if $\delta$ is a decreasing function such that $\delta(x)^3\log x$ tends to infinity with $x$ and $\cal P$ is $\delta$-dense, then by \rf{delta} and Theorem \ref{thm1} we see that  $K_N(\cal P)$ does not remain bounded. Hence, Theorem \ref{thm0} is merely a corollary
of Theorem \ref{thm1}.

We are not aware of earlier results of the type considered in Theorem \ref{thm0} or Theorem \ref{thm1} for sequences that are not quite dense. For other developments of the ideas deriving from \cite{ER,ET,R} see Pomerance \cite{Pom}.

With the sequence of primes comprising $\cal P$ we associate their {\em second differences}
\be{5} \Delta_n = p_{n+2} -2 p_{n+1} + p_n. \ee
Following R\'enyi in spirit, our approach to Theorem \ref{thm1} rests on  the observation that $\Delta_n$ is not too small for many values of $n$. Our next theorem is a strong quantitative version of this principle.

\begin{theorem}\label{thm2} Fix $x_0$ and $\delta$ as in Theorem \ref{thm1}. Then
there is
 a sequence of natural numbers $N_0(q)$ 
with the property that for all $N\ge N_0(q)$ and for all sets of primes  
$\cal P$ that are  $\delta$-dense relative to $x_0$ and some ${\cal P}_{q,a}$, 
one has
$$  \sum_{N<n\le 2N} \f{|\Delta_n|}{p_n} \le \f{11}{\delta_{2N+2}}. $$
If $\delta(x)^2\log x$ tends to infinity with $x$, then one also has
$$ \sum_{N<n\le 2N} \f{|\Delta_n|}{p_n} \ge 10^{-7} \delta_{2N}^3 . $$
\end{theorem}

Perhaps it is worth remarking that the upper bound
recorded in Theorem~\ref{thm2} is nearly the best possible. We demonstrate this with a scattered sequence that we briefly discuss at the end of the paper.

The proof of Theorem \ref{thm2} invokes upper bounds for the number of triplets of primes that come close 
to an arithmetic 3-progression. In Section 2 we use Selberg's sieve and a method of Gallagher \cite{G} to manufacture a suitable estimate, but it is worth pointing out that the older Brun's sieve and a technique of Hardy and Littlewood \cite{PN} would yield results of comparable strength. Equipped with the sieve estimates, the transition to the lower bound announced in Theorem \ref{thm2} is elementary, and is performed in Section 3.  For Theorem \ref{thm1}, we need a more explicit version of Theorem~\ref{thm2} (see Lemma \ref{L5} below) and the method of R\'enyi \cite{R}. The latter depends, in its original form, on the prime number theorem and is therefore not directly applicable to  subsets of the primes. In Section 4, we reconfigure R\'enyi's approach and establish Theorem \ref{thm1}. Thus our arguments that have elements in common with the work of Erd\H os and R\' enyi \cite{ER}, rely on methods that have been familiar for decades and yet, are of strength sufficient to address sequences of primes that are not quite dense. 

If more is known about the distribution of the sequence $(p_n)$, then our arguments
sometimes produce estimates that are superior to those recorded in the theorems.
For example this is the case when the $\delta$-dense set $\cal P$ has the additional
property that the numbers $p_{2n}/p_n$ remain bounded. With this extra assumption, the factor $\delta^{-1}$ can be deleted from the upper bounds in Theorems 2 and~3. 
For  more details on this refinement, the reader is referred to sections 3 and 4 below.

\section{A sieve estimate}

In this section we establish an auxiliary estimate concerned with triplets of primes. The main result is Lemma \ref{L2} below, and this depends on a certain singular series average that we now describe.

Throughout this section, let ${\mathbf h}=(h,h')\in\NN^2$ and suppose that $h'<h$. Then
\be{21} D=hh'(h-h') \ee
is an even  natural number. For a prime $p$, let $\nu_{\hh}(p)$ denote the number of distinct residue 
classes, modulo $p$, in which the numbers $0, h, h'$ lie. Then $1\le \nu_\hh(2)\le 2$ and
$1\le \nu_\hh(p) \le 3$ for all odd primes $p$. Further, it is immediate that one has 
$\nu_\hh(p) =3$ if and only if $p\nmid D$. For a given $q\in\NN$, we now define the 
{\em singular product}
\be{23} {\fk S}_{q,\hh} = \prod_{p\nmid q}\Big(1-\f{\nu_\hh(p)}p\Big)\Big(1-\f{1}p\Big)^{-3}. \ee
Note that
\be{24} 
       \Big(1-\f{\nu_\hh(p)}p\Big)\Big(1-\f{1}p\Big)^{-3} = 1+ a(p,\nu_\hh(p))
\ee
where
$$ a(p,\nu) = \f{p^3-(p-1)^3 - \nu p^2}{(p-1)^3}. $$
One readily checks that for odd primes one has
\be{stern} |a(p,\nu)| \le 3/(p-1) \quad (\nu=1,2), \quad |a(p,3)|\le 4/(p-1)^2.\ee
Hence, recalling that $\nu_\hh(p)=3$ holds for all $p\nmid D$, one finds that the product \rf{23} converges absolutely to a non-negative limit. 

\begin{lemma}\label{L1} Let $\eps>0$. Then,
uniformly for $q\in\NN$ and $0< \al\le 1$ one has
$$ \sum_{1\le h\le H} \sum_{|2h'-h|< \al h} {\fk S}_{q,\hh} =\f12  \al H^2 + O(H^{1+\eps}). $$
\end{lemma}

When $q=1$ a similar  estimate occurs in Gallagher \cite{G}, but the average there is over more parameters, and is more symmetric. We therefore give a complete proof, although we shall follow \cite{G} quite closely. 
For convenience, it is appropriate to put
$$ {\cal H} = \{(h,h')\in\NN^2: h\le H, \, |2h'-h|< \al h\}. $$
Note that $\al\le 1$ ensures that for any pair $(h,h')\in\cal H$ one
has $h'<h$. In particular, ${\fk S}_{q,\hh}$ is defined.
Now put
\be{aa} a_\hh(r) = \prod_{p|r} a(p,\nu_\hh(p)). \ee
Then, by \rf{24}, the absolutely convergent product in \rf{23} can be rewritten as
\be{25} {\fk S}_{q,\hh} = \multsum{r=1}{(r,q)=1}^\infty \mu(r)^2 a_\hh(r). \ee
From \rf{stern}, \rf{aa} and a familiar divisor estimate, we infer that $a_\hh(r) \ll r^{\eps-2} (r,D)$. Consequently, since $\hh\in\cal H$ implies  $0<D\le H^3$,  we conclude that
$$ \sum_{r>H} \mu(r)^2 |a_\hh(r)| \ll  H^{\eps-1} $$
holds uniformly in $\hh\in\cal H$. Then, the crude bound $\#{\cal H}\le H^2$ and \rf{25} suffice to deduce that 
\be{27} \sum_{\hh \in\cal H} {\fk S}_{q,\hh} = \multsum{r\le H}{(r,q)=1} \mu(r)^2 \sum_{\hh\in\cal H} a_\hh (r) +
O(H^{1+\eps}). \ee 
Note that this estimate is uniform with respect to $\alpha$ and $q$.

Consider the inner sum over $\hh$ in \rf{27} for a given
square-free number $r\le H$. Let $r=p_1\cdots p_\omega$ be the prime factorization.  We apply \rf{aa} and sort 
the $\hh\in\cal H$ according to given values of $\nu_\hh(p_j)$ $(1\le j\le \omega)$ to conclude that
\be{28} \sum_{\hh\in\cal H} a_\hh (r) = \multsum{1\le \nu_j\le 3}{1\le j\le \omega} a(p_1,\nu_1)\cdots a(p_\omega,\nu_\omega) S(r,\bfnu)\ee
where $S(r,\bfnu)$ is the number of $\hh\in\cal H$ with $\nu_\hh(p_j)=\nu_j$ for all $1\le j\le \omega$.
Note that the condition $\nu_\hh(p_j)=\nu_j$  depends only on the residue classes of $h$ and $h'$, modulo $p_j$. Hence, we may arrange $h$ and $h'$ into residue classes, modulo $r$, and then apply the Chinese Remainder Theorem to see that
$$ S(r,\bfnu) = \trisum{1\le a,a' \le r}{\nu_{\mathbf a}(p_j)=\nu_j}{1\le j\le \omega}
\trisum{\hh \in\cal H}{h\equiv a\bmod r}{h'\equiv a'\bmod r} 1. $$ 
For $r\le H$,  we also have   
$$ \trisum{\hh \in\cal H}{h\equiv a\bmod r}{h'\equiv a'\bmod r} 1
= \multsum{1\le h\le H}{h\equiv a\bmod r} \Big(\f{\al h}{r} +O(1)\Big) = \f12 \al \Big(\f{H}{r}\Big)^2 + 
O\Big(\f{H}{r}\Big). $$
Now let $t(p,\nu)$ denote the number of choices for $a, a'$ with $1\le a, a' \le p$ such that the numbers
$0,a,a'$ lie in exactly $\nu$ residue classes, modulo $p$. Then, again by the Chinese Remainder Theorem,
$$ \trisum{1\le a,a' \le r}{\nu_{\mathbf a}(p_j)=\nu_j}{1\le j\le \omega} 1 = \prod_{j=1}^\omega t(p_j,\nu_j),$$
and on collecting together we infer that
$$  S(r,\bfnu) = \Bigg(
\f12 \al \Big(\f{H}{r}\Big)^2 + 
O\Big(\f{H}{r}\Big) \Bigg) \prod_{j=1}^\omega t(p_j,\nu_j)      . $$
Now \rf{28} delivers
\be{29} \sum_{\hh\in\cal H} a_\hh (r) =  \f12 \al \Big(\f{H}{r}\Big)^2 \prod_{p\mid r} \sum_{\nu=1}^3 a(p,\nu)t(p,\nu) + O\Big(\f{H}{r}\prod_{p\mid r}\sum_{\nu=1}^3|a(p,\nu)|t(p,\nu) \Big). \ee

An inspection of the definition of $t(p,\nu)$ readily shows that $$t(p,1)=1, \quad t(p,2)= 3(p-1), \quad t(p,3)=(p-1)(p-2).$$ A short calculation leads to the identity
$$  \sum_{\nu=1}^3 a(p,\nu)t(p,\nu) = 0 $$
for all primes $p$, and for odd primes, by \rf{stern} we also have
$$  
\sum_{\nu=1}^3 |a(p,\nu)|t(p,\nu) \le 15. 
$$
It follows that the leading term in \rf{29} vanishes except when $r=1$. Moreover,
again using a  divisor estimate, we see  that
the error term  in \rf{29} does not exceed $O(Hr^{\eps-1})$. Hence, by \rf{27},
$$ \sum_{\hh\in\cal H} \fk{S}_{q,\hh} = \f12  \al H^2 + O(H^{1+\eps}) + O\Big(H \sum_{r\le H} r^{\eps-1}\Big), $$
and the conclusion of Lemma \ref{L1} follows.

\begin{lemma}\label{L2} Suppose that $0<\al\le 1\le H\le x$ and that $a,q\in\NN$
are coprime with $1\le a\le q$. Let $U=U_{\al,q,a}(x,H)$ denote the number of 
primes $p,p',p''$ with $p\equiv p'\equiv p''\equiv a\bmod q$ that satisfy 
the inequalities
\be{pcond} 5\le p\le x,\quad p< p'' \le p+qH, \quad |p''-2p'+p| < \al (p''-p). \ee
Further let $\eps>0$. 
Then there are a number $x_2=x_2(q)$ depending only on $q$ 
and a number $E=E_\eps$ depending only on $\eps$ such that whenever $x\ge x_2$ one has  
\be{triple} U \le (25\al H^2 +  E H^{1+\eps})) \f{x}{\phi(q)(\log x)^3}. \ee
\end{lemma} 

{\em Proof}. Suppose that $p,p',p''$ is a triple counted by $U$. We write
\be{X0} p= a+ql,\quad p''-p=qh, \quad p'-p = qh'. \ee
Then $l\in\NN_0$, $(h,h')\in\NN^2$, and the conditions \rf{pcond} imply that
\be{X1} 0\le l\le x/q, \quad h\le H, \quad |h-2h'|\le \al h. \ee
By \rf{X0}, it follows that $U$ does not exceed the number of 
 $l\in\NN_0$, $(h,h')\in\NN^2$ satisfying \rf{X1} and $a+ql\ge 5$ for which the three numbers
\be{X2} 
a+ql,\quad a+q(l+h), \quad a+q(l+h') \ee
are all prime. 

Let $V(h,h')=V(\hh)$ denote the number of integers $l$ with $0\le l\le x/q$ and $a+ql\ge 5$
for which the numbers \rf{X2} are simultaneously prime. Then, in the notation
of the proof of Lemma \ref{L1}, the above argument shows that
$$ U\le \sum_{\hh\in\cal H} V(\hh).$$
Further, the quantity $V(\hh) $ is readily estimated by an upper bound sieve. We wish to apply  \cite[Theorem 5.7]{HR}, and with this end in view we consider, for a prime $p$, the number $\rho_\hh(p)$ of incongruent solutions in $z$ of the congruence
$$ (a+qz)(a+q(z+h))(a+q(z+h')) \equiv 0 \bmod p.$$
Then, whenever $p\mid q$, one has $\rho_\hh(p)=0$ while in the contrary case
$p\nmid q$ it is immediate that $\rho_\hh(p) = \nu_\hh (p)$. If $\hh\in\cal H$ is such that $\rho_\hh(p)<p$
holds for all primes $p$, then 
\cite[Theorem 5.7]{HR} is applicable and delivers the inequality 
\be{X5} V(\hh) \le 50 {\fk S}_{q,\hh} \f{x}{\phi(q)}(\log x)^{-3}\ee
for all $x$ that are sufficiently large in terms of $q$, as one readily confirms by inspecting
\rf{23} and the Euler product in \cite[(5.8.3)]{HR}. 

It remains to evaluate $V(\hh)$ in those cases where $\rho_\hh(p)<p$ fails for some prime $p$.
The trivial upper bound $\rho_\hh(p)\le \min (3,p)$ shows that this is possible only when $p=2$ or $3$.
Further, the hypothesis that $\rho_\hh(2)=2$ implies that $2\nmid q$, and that at least one of $h,h'$
is odd. By \rf{X1} we then find that one of the differences $p'-p$, $p''-p$ is odd which is impossible for $p\ge 5$.
This shows that  $\rho_\hh(2)=2$ implies $V(\hh)=0$, and a similar argument confirms that the same is 
true when $\rho_\hh(3)=3$. In particular, we now see that \rf{X5} holds for all $\hh\in\cal H$. Summing \rf{X5}
over these $\hh$ with the aid of Lemma \ref{L1} yields Lemma \ref{L2}.

\section{Second differences - Proof of Theorem 3}

We launch an attack toward the estimates claimed in Theorem \ref{thm2} with a preliminary remark. 
Throughout, suppose that $x_0$ and $\delta$ are fixed, as in Theorem \ref{thm1}. Let $\cal P$ be
a set of primes, choose $a,q$ with ${\cal P}\subset{\cal P}_{q,a}$, and assume that \rf{3} holds for all
$x\ge x_0$. 
Suppose it were the case that $p_n>n^2$ holds for some $n$ with $n^2\ge x_0$. Then, in \rf{3}
we
take $x=n^2$ and use the lower bound for $\delta(x)$ to infer that
$$ n \ge 
\delta(n^2)\pi(n^2;q,a) \ge (2\log n)^{-1}\pi(n^2;q,a). $$
The prime number theorem in arithmetic progressions supplies a number $x_1(q)$ such that whenever
$x\ge x_1(q)$ then one has $\pi(x;q,a)\ge x/(2\phi(q)\log x)$. Hence, for $n^2\ge \max(x_0,x_1(q))$, we conclude that 
$$ n \ge \f{n^2}{4\phi(q)(\log n)^2}.  $$
This is absurd for $n$ sufficiently large in terms of $q$. 
It follows that there is a number $n_0$, depending only on $x_0$ and $q$, with the property that whenever $n\ge n_0$ then the inequalities
\be{crude} p_n\le n^2 \quad \mbox{and} \quad \delta_n^{-1} \le 2 \log n \ee
hold.
These bounds are improved in the following lemma, but they play a role in its proof.

\begin{lemma}\label{L4} Let $x_0, \delta, {\mathcal P}$ and $a,q$ be as in the preceding paragraph.
Then there is a number $n_0$ depending only on $x_0$ and $q$ such that whenever $n\ge n_0$, one has
$$ \f34 \phi(q)n\log n \le p_n \le 2 \phi(q)\delta_n^{-1} n \log n. $$ 
\end{lemma}

Within the proof, we may suppose that \rf{3} holds with $\delta=\delta(x)$. 
But then, for $x\le p_n$, the bound \rf{3} also holds with $\delta=\delta_n$.
Now suppose for contradiction that
$p_n>x_0$ and 
$p_n> 2\phi(q)\delta_n^{-1} n\log n$ hold simultaneously. We may use \rf{3} with $x=  2\delta_n^{-1}\phi(q) n\log n$ and  
then see  that
$$  n \ge  \pi( 2\delta_n^{-1}\phi(q) n\log n;q,a).$$
Using the prime number theorem in arithmetic progressions much as above, this implies via \rf{crude} that
$$n\ge \f32  n \f{\log n}{\log \phi(q)n}.$$ This is certainly false for $n$ large in terms of $q$.
The upper bound for $p_n$ follows.

Next, let $\varpi$ denote the $n$-th member of the ascending sequence of {\em all} primes in ${\cal P}_{q,a}$. 
Then $p_n\ge \varpi$,
and by the prime number theorem in arithmetic progressions once again, 
one has $\varpi\ge \f34\phi(q)n\log n$ for all large $n$. This
completes the proof of Lemma \ref{L4}. 

\medskip

The lower bound \rf{delta} is now immediate. Indeed, by Lemma \ref{L4} and \rf{crude},
we have
$$ \delta_n = \delta(p_n) \ge \delta(2\phi(q)\delta_n^{-1} n\log n) \ge \delta (4\phi(q) n (\log n)^2), $$
as required.

\medskip

The next task ahead of us is to establish Theorem \ref{thm2}. For the upper bound, we apply the triangle inequality to \rf{5} and then see from Lemma \ref{L4} that
\be{tele} \sum_{N<n\le 2N} \f{|\Delta_n|}{p_n} \le 2 \sum_{N<n\le 2N+1}
\f{p_{n+1}-p_n}{p_{N+1}} \le \f{2p_{2N+2}}{p_{N+1}} \le  \f{11}{\delta_{2N+2}}, \ee 
provided only that $N$ is large. This already completes the proof of the upper bound, but there is a simple variant of this argument. Suppose the $\delta$-dense set has the additional property that there exist a number $A$ such that 
\be{bounded} \f{p_{2m}}{p_m} \le A \ee
holds for all large $m$. Then the inequalities in \rf{tele} provide the alternative estimate
\be{tele2} \sum_{N<n\le 2N} \f{|\Delta_n|}{p_n} \le 2A. \ee
This substantiates a remark that we have made in the introductory part of the paper.

\medskip
The lower bound in Theorem \ref{thm2} will be deduced from the sieve bounds established in the previous section.
It will be useful to introduce the parameters
\be{para}  C= 2/\delta_{2N}, \quad B= 10^{-5} C^{-2}. \ee
We will use this notation in the remainder of this paper.

\begin{lemma}\label{L5} 
Fix $x_0$ and $\delta$ as in Theorem \ref{thm1}, and suppose that
$\delta(x)^2\log x$ tends to infinity with $x$. Then
there is
 a sequence of natural numbers $N_2(q)$ 
with the property that  for all sets of primes  
$\cal P$ that are  $\delta$-dense relative to $x_0$ and some ${\cal P}_{q,a}$, 
and for all $N\ge N_2(q)$, the set
$$ {\cal B}(N)=\{N<n\le 2N-2:\, p_{n+2}-p_n\le 33C\phi(q)\log N,\, |\Delta_n|\ge B\phi(q) \log N\} $$
has at least $N/2$ elements.
\end{lemma}

Note that once this lemma is established, we may apply Lemma \ref{L4} to conclude that
$$ \sum_{n\in  {\cal B}(N)} \f{|\Delta_n|}{p_n}\ge  \f14 C^{-1}B > 10^{-7} \delta_{2N}^3 $$
holds whenever $N\ge N_2(q)$. This includes  
the lower bound recorded in Theorem \ref{thm2}.  

\medskip 

 We now give the  proof of Lemma \ref{L5}.
Let ${\cal S}_1$ be the set of all $n\in(N,2N-2]$ where $p_{n+2}-p_n > 33 C \phi(q) \log N$. Then 
$$ \# {\cal S}_1 \le \sum_{N<n\le 2N} \f{p_{n+2}-p_n}{ 33 C \phi(q) \log N} \le \f{p_{2N}+p_{2N-1}}{ 33 C \phi(q) \log N},$$
and by Lemma \ref{L4}, for sufficiently large $N$, we conclude that
$ \# {\cal S}_1 \le N/8$.

Next, 
let ${\cal S}_2$ be the set of all $n\in(N,2N-2]$ where $p_{n+2}-p_n \le B \phi(q) \log N$. Then, in view of
\rf{5} and
Lemma \ref{L4}, the primes $p=p_n$, $p'=p_{n+1}$ and $p''=p_{n+2}$ satisfy the conditions \rf{pcond} with $\al=1$, $x=2C \phi(q)N \log 2N$ and $H=\f{\phi(q)}{q}B\log N$. Note that \rf{delta}  implies that $B\log N$ tends to infinity with $N$.
Hence  $H$ is large when $N$ is large (in terms of $q$). Therefore, for large $N$,
Lemma \ref{L2} is applicable, and delivers the estimates  
\begin{align*}
\# {\cal S}_2 &\le U_{1,q,a}(2C \phi(q)N \log 2N, H) \\
& \le (25H^2+E_\eps H^{1+\eps}))\f{2C \phi(q)N \log 2N}{( \phi(q)\log 2N + O(\log \log N))^3} \\
& \le 100 B^2 C N \le 10^{-8}N,
\end{align*}
provided again that $N$ is sufficiently large in terms of $q$. 
 
For notational convenience,  put 
$$ J_0= B\phi(q)\log N, \quad J_1 = 33 C\phi(q) \log N, $$
and let  ${\cal S}_3$ be the set of all $n\in(N,2N-2]$ where 
$$ J_0 < p_{n+2}-p_n \le J_1 \quad \mbox{and}\quad |\Delta_n|<J_0. $$
Note that a number $n\in(N,2N-2]$
that is in none of the sets ${\cal S}_1$, ${\cal S}_2$, ${\cal S}_3$ lies in ${\cal B}(N)$.
Hence, once we have proved that $\# {\cal S}_3\le N/10$ holds for large $N$, the proof of 
Lemma \ref{L5} will be complete.  

We proceed by a dissection argument. Let ${\cal S}_3(H)$ be the subset of ${\cal S}_3$ where $\f12 qH<p_{n+2}-p_n\le qH$.
The range of $H$ relevant to us is the interval
$$ 2J_0 \le qH \le 2J_1. $$
For such a value of $H$ and $n\in{\cal S}_3(H)$ we have $|\Delta_n|\le J_0 = \al (\f12qH) \le \al (p_{n+2}-p_n)$
where we chose
$$ \al = \f{2J_0}{qH}. $$ 
By Lemma  \ref{L4}, we may take 
$x= 2C\phi(q) N\log 2N$ in Lemma 
 \ref{L2}, and then find that
\begin{align*}
\# {\cal S}_3(H) &\le U_{\al,q,a}(2C\phi(q) N \log 2N, H) \\
&\le (25 \al H^2+ E_\eps H^{1+\eps}) \f{2CN\log 2N}{(\log(2C\phi(q) N\log 2N))^3}. 
\end{align*}
When $N$ is sufficiently large in terms of $q$, this upper bound simplifies to
$$ \# {\cal S}_3(H) 
\le (50  J_0Hq^{-1}+ E_\eps H^{1+\eps}) \f{3CN}{(\log 2N)^2} .$$
We take $H=2^j J_0q^{-1}$ with $2^j \le J_1/J_0 = 33C/B$ and sum over $j\ge 1$. Since  ${\cal S}_3$
is covered by the union of these ${\cal S}_3(H)$, we conclude that for large $N$ we have
$$
\# {\cal S}_3 
 \le (50J_0 J_1q^{-2}+ E_\eps(J_1/q)^{1+\eps}))\f{6CN}{(\log 2N)^2} \le 10^4 BC^2 N .
$$
By hypothesis, we have $10^4 BC^2 = 1/10$ so that indeed $\# {\cal S}_3\le N/10$, as required.

\section{Curvature - Proof of Theorem 2}
This section is devoted to the proof of Theorem \ref{thm1}. Recall that $z_n=n+{\mathrm i} \log p_n$, and with
this choice put
$$ k_n =  \Big| \arg \f{z_{n+1}-z_n}{z_n-z_{n-1}}\Big|.$$
By the elementary properties of the arctan function and its addition theorem, one finds the alternative representations
\begin{align} k_n &= \Big| \arctan \log \f{p_{n+2}}{p_{n+1}} - \arctan\log \f{p_{n+1}}{p_{n}}\Big|\notag \\ &
 = \arctan \frac{\displaystyle \Big| \log \f{p_{n+2}p_n}{p_{n+1}^2} \Big|}{\displaystyle 1+ \Big( \log \f{p_{n+2}}{p_{n+1}}\Big)\Big(\log \f{p_{n+1}}{p_{n}}\Big)}. \label{kk}
\end{align}
It will also be useful to put
\be{Gam} \Gamma_n = (p_{n+2}-p_{n+1}) (p_{n+1}-p_{n})p_{n+1}^{-2}. \ee
Then, by \rf{5}, one has
\be{31}  \f{p_{n+2}p_n}{p_{n+1}^2} = 1 + \f{\Delta_n}{p_{n+1}} - \Gamma_n. \ee
This identity also occurs in R\'enyi's work (see \cite{R}, eqn.\ (5)), and is crucial to his arguments.
Before we proceed further, we note that 
$$0 < \Gamma_n  \le  \f{p_{n+2}-p_{n+1}}{p_{n+1}} .  $$ In particular, Lemma \ref{L4} yields
$\Gamma_n \le    (p_{n+2}-p_{n+1})/ (\f34\phi(q)n\log n)$, and  the argument that was used in \rf{tele} produces
\be{Gamma} \sum_{N<n\le 2N} \Gamma_n \le \f{6}{\delta_{2N+2}},  \ee
at least when $N$ is large in terms of $q$, as we assume from now on. For
such $N$
we  proceed to derive the inequality
\be{upper} \sum_{N<n\le 2N} k_n \le\f{300}{\delta_{2N+2}} . \ee
Let $\cal C$ denote the set of all $n\in(N,2N]$ where $\Gamma_n>1/8$, and let 
$\cal D$ denote the set of all $n\in(N,2N]$ where $|\Delta_n|>p_n/8$. Finally, let $\cal E$ denote the
set of the remaining $n$ in $(N,2N]$. By \rf{Gamma}, it follows that $\#{\cal C} \le 48/ \delta_{2N+2} $, and similarly,
one finds from \rf{tele} that  $\#{\cal D} \le  88/\delta_{2N+2}$. The trivial bound $k_n\le \pi/2$ suffices to see
that the contribution from all $n\in {\cal C}\cup\cal D$ to the sum in \rf{upper} does not exceed $68\pi/\delta_{2N+2}$.

Now suppose that $n\in\cal E$. Then  
$$ \Big| \f{\Delta_n}{p_{n+1}} - \Gamma_n\Big| \le \f14. $$
Hence, on applying the familiar bound $\arctan t \le t$ that is valid for all $t\ge 0$, we first see
from \rf{kk} that
$$ k_n \le \Big| \log \f{p_{n+2}p_n}{p_{n+1}^2} \Big|, $$
and then, by observing that for real numbers $t$ with $|t|\le 1/4$ one has $|\log (1+ t)| \le 2|t|$, we conclude
via \rf{31} that
$$ k_n \le 2 \f{|\Delta_n|}{p_{n}} + 2\Gamma_n. $$
We sum this over $n\in\cal E$. Then, by \rf{tele} and \rf{Gamma}, we see that  \rf{upper} indeed holds.

\medskip
The upper bound reported in Theorem \ref{thm1} is readily deduced from \rf{upper}. There is a number
$N_1(q)$ such that \rf{upper} holds for all $N\ge N_1(q)$. But then, if a large $M$ is given, we may
take $N=2^jN_1(q)\le M$ in \rf{upper} and sum over $j$.
Using the trivial bound
for $k_n\le \pi/2$ when $n\le N_1(q)$ and observing that 
 $\delta$ is decreasing, we find in this way that
$$ K_M({\cal P}) \le\f{\pi}{2} N_1(q) +499 \f{\log M}{\delta_M}. $$
When $N$ is sufficiently large in terms of $q$, this implies the upper bound
recorded in Theorem \ref{thm1}. 
Again, there is a variant of this argument in the case where
\rf{bounded} holds. Then we have \rf{tele2} available, and for the same reason
in \rf{Gamma} the upper bound can be replaced by $A$.  With these estimates in hand, the above argument produces the better bound
$$ K_N({\cal P}) \ll_A \log N. $$
Once again, this confirms a claim from the introduction. 

\medskip

The verification of the  lower bound in Theorem \ref{thm1} is somewhat more complex. 
Throughout the argument below  we use the notation as introduced in   
Lemma \ref{L5}. Let $N\ge N_0(q)$, and consider a number
$n\in{\cal B}(N)$. Then by  \rf{5} and the triangle inequality,
$$ B\phi(q) \log N \le |\Delta_n| \le 33C\phi(q)\log N.$$ 
Furthermore, one has $p_{n+2}-p_{n+1}\le p_{n+2}-p_n\le 33C\phi(q)\log N$, and the same inequality holds for 
$p_{n+1}- p_n$. Hence, by Lemma \ref{L4} and \rf{Gam}, we see that $ \Gamma_n\le (44C)^2 N^{-2}$, and that
$$ \f{B}{2CN}\f{ |\Delta_n|}{p_{n+1}} \le \f{44C}{N}.$$
Here we have used that $N$ is large. These last inequalities combine with the bound on $\Gamma_n$ to
$$
\f{B}{3CN} \le \Big|  \f{\Delta_n}{p_{n+1}} - \Gamma_n \Big| \le   \f{45C}{N} \le \f19. $$
However, when $|t|\le \f12$ one has $\f12|t|\le |\log (1+t)| \le 2|t|$,
so that \rf{31} now yields 
$$
\f{B}{6CN} \le \Big| \log \f{p_np_{n+2}}{p^2_{n+1}} \Big| \le   \f{90C}{N} , $$

We apply Lemma \ref{L4} again to confirm that for $j=1$ and $2$, one has 
$$ 1\le \f{p_{n+j+1}}{p_{n+j}} = 1+ \f{p_{n+j+1}-p_{n+j}}{p_{n+j}}\le 1+ \f{44C}{N} \le \f{10}9. $$
Thus 
$$ 0\le \log  \f{p_{n+j+1}}{p_{n+j}} \le \f19, $$
and hence,
$$
\f{B}{7CN} \le  \frac{\displaystyle \Big| \log \f{p_{n+2}p_n}{p_{n+1}^2} \Big|}{\displaystyle 1+ \Big( \log \f{p_{n+2}}{p_{n+1}}\Big)\Big(\log \f{p_{n+1}}{p_{n}}\Big)} 
\le   \f{90C}{N} \le \f29. $$
For $0\le t\le \f14$ one has $\arctan t \ge \f12 t$. Therefore, by \rf{kk}, we conclude that
$k_n\ge B/(14CN)$ holds for all $n\in {\cal B}(N)$. By \rf{para} and Lemma \ref{L5} it follows that
$$ \sum_{n\in{\cal B}(N)} k_n \ge \f{B}{28 C} \ge\f1{10^8 \delta_{2N}^{3}}. $$
In this estimate, we replace $N$ by $2^{-j}N$ and sum over $1\le j\le \f13 \log N$. 
But then $2^{-j}N \ge \sqrt N$ for all $j$, and we only have to arrange that $\sqrt N \ge N_2(q)$,
with $N_2$ as in Lemma \ref{L5}. We conclude that
$K_N({\cal P}) \ge 10^{-8} \delta_N^{-3}\log N$, as required to complete the proof of Theorem \ref{thm1}. 

\section{A scattered sequence} 
We end with a brief description of a sequence with large curvature. Let $A_1=10$ and define $A_l$ by the recursion $A_{l+1} = 2 A_l \log 4A_l$ $(l\ge 1)$. Note that the intervals $I_l=[A_l, 4A_l]$ are disjoint. We now construct a set of primes $\cal Q$ as follows. If the number of primes in $I_l$ is even, then all these primes become elements of $\cal Q$, and in the contrary case, we put all but the smallest of
the primes in $I_l$ in $\cal Q$. Primes that are not in some $I_l$ are not in $\cal Q$.
Note that we have arranged that the number of elements in ${\cal Q}\cap I_l$ is even.

We claim that $\cal Q$ is $\delta$-dense with $\delta(x)=1/\log x$. To see this,
let $\pi_{\cal Q}(x)$ denote the number of primes in $\cal Q$ not exceeding $x$,
and suppose that $x$ is large. Then, there is some $l$ with $\f54 A_l \le x \le \f54 A_{l+1}$. In the case where $\f54 A_l \le x \le 4A_l$, Chebyshev's estimates give $\pi_{\cal Q}(x) \gg x/ \log x$ which is more than is required. In the range
$4A_l\le x \le \f54 A_{l+1}$ we use the prime number theorem to see that
\begin{align*} \pi_{\cal Q} (x) & \ge \pi_{\cal Q}(4A_l) \ge \f{3A_l}{\log 4A_l} (1+o(1)) \\
&= \Big(\f32 +o(1)\Big) \f{A_{l+1}}{(\log A_{l+1})^2} 
=    \Big(\f65 +o(1)\Big) \f{\f54 A_{l+1}}{(\log \f54 A_{l+1})^2}.
\end{align*}
In particular, this confirms \rf{3} with $\delta(x)=1/\log x$, as desired.

Let $(q_j)$ denote the sequence of the elements of $\cal Q$ in ascending order.
Now let $l$ be large. By construction, $\pi_{\cal Q}(4A_l)$ is even, say $2N$.
Then $q_{2N}<4A_l$ but $q_{2N+1}>A_{l+1}$. Also, by the prime number theorem, 
$q_{2N+2}-q_{2N+1}= o(A_{l+1})$ so that we now have $|\Delta_{2N}|\ge A_{l+1}(1+o(1))$,
and hence, again by the prime number theorem,
$$ \f{|\Delta_{2N}|}{q_{2N}} \ge \f{A_{l+1}}{4A_l}(1+o(1) = \Big(\f12+o(1)\Big)
\log A_l. $$
Further, the equation $2N=\pi_{\cal Q}(4A_l)$ and the straightforward  bounds
$$ \f{3A_l}{\log A_l}(1+o(1)) \le \pi_{\cal Q}(4A_l) \le  \f{4A_l}{\log A_l}(1+o(1)) $$
imply that $\log A_l = (1+o(1)) \log N$, so that we arrive at
$$ \f{|\Delta_{2N}|}{q_{2N}} \ge \f13 \log N. $$
In particular, we see that the sum considered in Theorem \ref{thm2} contains a single term exceeding $\f13 \log N$, which is of the order of $\delta_{2N+2}^{-1}$. 


\begin{thebibliography}{10}
\bibitem{ER} Erd\H os, P.; R\' enyi, A. Some problems and results on consecutive primes. Simon Stevin 27 (1950), 115--125.
\bibitem{ET} Erd\H os, P.; Tur\' an, P. On some new questions on the distribution of prime numbers. 
Bull. Amer. Math. Soc. 54 (1948), 371--378.
\bibitem{G} Gallagher, P. X. On the distribution of primes in short intervals. Mathematika 23 (1976), 4--9.
\bibitem{HR} Halberstam, H.; Richert, H.-E. Sieve methods. London Mathematical Society Monographs, No. 4. Academic Press, London-New York, 1974.
\bibitem{PN} Hardy, G.H.; Littlewood, J.E. Some problems of ???Partitio Numerorum???: III. On the expression of a number as a sum of primes, Acta Math. 44 (1922), 
1--70. 
\bibitem{Pom} Pomerance, C. The prime number graph. Math. Comp. 33 (1979), 399--408.
\bibitem{P} Prachar, K.
Bemerkung zu einer Arbeit von Erd\H os und R\' enyi und Berichtigung. 
Monatsh. Math. 58 (1954), 117.
\bibitem{R} R\'enyi, A. On a theorem of Erd\H os and Tur\' an. Proc. Amer. Math. Soc. 1 (1950), 7--10. 



\end{thebibliography}
\end{document}